\documentclass{article}

\usepackage{amsmath,amssymb,amsthm}
\usepackage[T1]{fontenc}
\usepackage{latexsym}

\def \vs{\vskip .6cm}

\def \beq{\begin{eqnarray*}}
\def \eeq{\end{eqnarray*}}

\def \nb{\overline \nabla}

\def \Rb{\overline R}

\def \h{\mathfrak h}
\def \g{\mathfrak g}

\def \m{\mathfrak m}

\def \RM{\mathbb{R}}

\def \ZM{\mathbb{Z}}
\def \CM{\mathbb{C}}
\def \HM{\mathbb{H}}

\def \LLambda {{\bf \Lambda}}
\def \llambda {\boldsymbol{\lambda}}

\def \leftr {[\hbox{\hspace{-0.15em}}[}
\def \rightr {]\hbox{\hspace{-0.15em}}]}
\def \la {\langle}
\def \ra {\rangle}

\newtheorem{defi}{Definition}[section]

\newtheorem{prop}[defi]{Proposition}

\newtheorem{conj}{Conjecture}

\newtheorem{theo}{Theorem}

\newtheorem{theo'}[defi]{Theorem}

\newtheorem{lemm}[defi]{Lemma}

\newtheorem*{NB}{Remark}

\newtheorem{coro}[defi]{Corollary}

\title{Homogeneous nearly Kähler manifolds}

\date{}

\author{Jean-Baptiste Butruille}

\begin{document}

\maketitle

\begin{abstract} 
We classify six-dimensional homogeneous nearly Kähler manifolds and give a positive answer to Gray and Wolf's conjecture: every homogeneous nearly Kähler manifold is a Riemannian 3-symmetric space equipped with its canonical almost Hermitian structure. The only four examples in dimension 6 are $S^3 \times S^3$, the complex projective space $\CM P^3$, the flag manifold $\mathbb F^3$ and the sphere $S^6$. We develop, about each of these spaces, a distinct aspect of nearly Kähler geometry and make in the same time a sharp description of its specific homogeneous structure.
\end{abstract}

%\ni {\it Keywords: } 

%\ni {\it AMS Classification: }

\section*{Introduction}

Probably the first example known of a nearly Kähler manifold is the round sphere in dimension 6, equipped with its well-known non integrable almost complex structure, introduced in \cite{fu}. The resulting \emph{almost Hermitian structure} is invariant for the action of $G_2$ on $S^6$ coming from the octonions (we look at $S^6$ as the unit sphere in the imaginary set $\Im \subset \mathbb O$). Thus $S^6 \simeq G_2/SU(3)$ is an example of a \emph{6-dimensional homogeneous nearly Kähler manifold} as we consider in the present article. Notice also that the representation of $SU(3)$ on the tangent spaces -- the isotropy representation -- is irreducible. Then, Joseph A. Wolf in his book {\it Spaces of constant curvature} \cite{wo}, discovered a class of isotropy irreducible homogeneous spaces $G/H$ that generalizes $S^6$. Indeed, $G$ is a compact Lie group and $H$, a maximal connected subgroup, centralizing an element of order 3. Further in this way, Wolf and Gray \cite{wo} investigated the homogeneous spaces defined by Lie group automorphisms i.e. such that $H$ is the fixed point set of some $s : G \to G$. They asked the following question : which of these spaces admit an invariant almost Hermitian structure with good properties ? In particular, nearly Kähler manifolds are associated, through their work, to a type of homogeneous spaces -- the 3-symmetric spaces -- corresponding to $s$ of order 3. Since their study was quite general, they felt legitimate to make a conjecture, which I reformulated using the terminology of the later article \cite{gr2} by Gray alone:

\begin{conj}[Gray and Wolf]
Every nearly Kähler homogeneous manifold is a 3-symmetric space equipped with its canonical almost complex structure.
\label{conjecture}
\end{conj}

Another way to construct examples is by twistor theory. The twistor space $Z$ of a self-dual 4-manifold has a natural complex structure by \cite{at} but also, a non integrable almost complex structure (see \cite{ee}). The latter may be completed into a nearly Kähler structure when the base is furthermore Einstein, with positive scalar curvature. The same construction holds for the twistor spaces of the positive quaternion-Kähler manifolds (see \cite{al,na}) or of certain symmetric spaces (as explained for instance in \cite{sa2}). In this case $Z$ is a 3-symmetric space again.

Next, in the 70's \cite{gr,gr3} and more recently (we mention, as a very incomplete list of references on the topic: \cite{re,be,na} and \cite{na2}), nearly Kähler manifolds have been studied for themselves. Some very interesting properties were discovered, especially in dimension 6, that give them a central role in the study of special geometries with torsion. These properties can all be interpreted in the setting of \emph{weak holonomy}. Nearly Kähler 6-dimensional manifolds are otherwise called \emph{weak holonomy $SU(3)$}. Many definitions of this notion have been proposed. One is by spinors \cite{ba} ; others by differential forms \cite{se,bu3} ; the original one by Gray \cite{gr5} was found incorrect (see \cite{al2}). Finally, Cleyton and Swann \cite{cl,cl2} have explored theirs based on the torsion of the canonical connection of a $G$-structure. They lead to a theorem which can be applied to our problem to show that nearly Kähler manifolds whose canonical connection has \emph{irreducible} holonomy are either 3-symmetric or 6-dimensional. This constitutes an advance towards and gives a new reason to believe to Gray and Wolf's conjecture. However, dimension 6 still resists. Moreover, when the holonomy is reducible, we do not have of a de Rham-like theorem, like in the torsion-free situation. It was Nagy's main contribution to this issue, to show that we can always lead back, in this case, to the twistor situation. Using this, he was able to reduce conjecture \ref{conjecture} to dimension 6. 

This is where we resume his work. We classify 6-dimensional nearly Kähler homogeneous spaces and show that they are all 3-symmetric.
\begin{theo}
Nearly Kähler, 6-dimensional, Riemannian homogeneous manifolds are isomorphic to a finite quotient of $G/H$ where the groups $G$, $H$ are given in the list:
\begin{itemize}
\item[--] $G=SU(2) \times SU(2)$ and $H=\{1\}$
\item[--] $G=G_2$ and $H=SU(3)$. In this case $G/H$ is the round 6-sphere.
\item[--] $G=Sp(2)$ and $H=SU(2)U(1)$. Then, $G/H \simeq \CM P^3$, the 3-dimensional complex projective space.
\item[--] $G=SU(3)$, $H=U(1) \times U(1)$ and $G/H$ is the space of flags of $\CM^3$.
\end{itemize}
Each of these spaces, $S^3 \times S^3$, $S^6$, $\CM P^3$ and $\mathbb F^3$, carry a unique invariant nearly Kähler structure, up to homothety.
\label{classification}
\end{theo}

As a corollary,

\begin{theo}
Conjecture \ref{conjecture} is true in dimension 6 and thus, by the work of Nagy, in all (even) dimensions.
\label{resolution conjecture}
\end{theo}

The proof is systematic. We start (proposition \ref{algebraic classification}) by making a list of the pairs $(G,H)$ such that $G/H$ is likely, for topological reasons, to admit an invariant nearly Kähler structure. The homogeneous spaces that appear in this list are the four compact 6-dimensional examples of 3-symmetric spaces found in \cite{wo}. Then, we show, for each space, that there exists no other homogeneous nearly Kähler structure than the canonical almost complex structure on it. 

As a consequence, this article is mainly focused on examples. However, the four spaces in question are quite representative of a number of features in nearly Kähler geometry. Thus, it might be read as a sort of survey on the topic (though very different from the one written by Nagy in the same handbook). Section 2 contains the most difficult point in the proof of theorem \ref{classification}. We look for left-invariant nearly Kähler structures on $S^3 \times S^3$. For this we had to use the algebra of differential forms on the manifold, nearly Kähler manifolds in dimension 6 being characterized, by the work of Reyes-Carri\'on \cite{re} or Hitchin \cite{hi3}, by a differential system on the canonical $SU(3)$-structure. Section 3 is devoted to two homogeneous spaces $\CM P^3$ and $\mathbb F^3$ that are the twistor spaces of two 4-dimensional manifolds: respectively $S^4$ and $\CM P^2$. We take the opportunity to specify the relation between nearly Kähler geometry and twistor theory. Weak holonomy stands in the background of section 4 on the 6-sphere. Indeed, we may derive this notion from that of \emph{special holonomy} trough the construction of the Riemannian cone. Six-dimensional nearly Kähler -- or weak holonomy $SU(3)$ -- structures on $S^6$ are in one-to-one correspondance with constant 3-forms, inducing a reduction of the holonomy to $G_2$ on $\RM^7$. Finally, in section 1, we introduce the notions and speak of 3-symmetric spaces in general and in section 5, we provide the missing elements for the proof of theorems \ref{classification} and \ref{resolution conjecture}.

We should mention here the remaining conjecture on nearly Kähler manifolds:
\begin{conj}
Every \emph{compact} nearly Kähler manifold is a 3-symmetric space.
\end{conj}
That this conjecture is still open means in particular that the homogeneous spaces presented along this article are the only known compact (or equivalently, complete) examples in dimension 6. Again by the work of Nagy \cite{na2}, it may be separated in two restricted conjectures. The first one relates to a similar conjecture on quaternion-Kähler manifolds and symmetric spaces, for which there are many reasons to believe that it is true (to begin with, it was solved by Poon and Salamon \cite{po} in dimension 8 and recently by Haydeé and Rafael Herrera \cite{herr} in dimension 12). The second may be formulated: the only compact, simply connected, \emph{irreducible} (with respect to the holonomy of the intrinsic connection), \emph{6-dimensional}, nearly Kähler manifold is the sphere $S^6$ -- and concerns the core of the nearly Kähler geometry: the fundamental explanation of the rareness of such manifolds or the difficulty to produce non-homogeneous examples.

\section[Preliminaries]{Preliminaries: nearly Kähler manifolds and 3-symmetric spaces}

Nearly Kähler manifolds are a type of \emph{almost Hermitian manifolds} i.e. $2n$-dimensional real manifolds with a $U(n)$-structure (a $U(n)$-reduction of the frame bundle) or equivalently, with a pair of tensors $(g,J)$ or $(g,\omega)$, where $g$ is a Riemannian metric, $J$ an almost complex structure compatible with $g$ in the sense that 
\[ \forall X,Y \in TM, \quad g(JX,JY)=g(X,Y) \]
($J$ is orthogonal with respect to $g$ pointwise) and $\omega$ is a differential 2-form, called the \emph{Kähler form}, related to $g,J$ by
\[ \forall X,Y \in TM, \quad \omega(X,Y)=g(JX,Y) \]

Associated with $g$ there is the well-known Levi-Civita connection, $\nabla$, metric preserving and torsion-free. But nearly Kähler manifolds, as every almost Hermitian manifolds, have another natural connection $\nb$, called the \emph{intrinsic connection} or the \emph{canonical Hermitian connection}, which shall be of considerable importance in the sequel. Let $\mathfrak{so}(M)$ be the bundle of skew-symmetric endomorphisms of the tangent spaces (the adjoint bundle of the metric structure). The set of metric connections of $(M,g)$ is an affine space  $\mathcal{SO}$ modelled on the space of sections of $T^*M \otimes \mathfrak{so}(M)$. Then, the set $\mathcal U$ of Hermitian connections (i.e. connections which preserve both the metric and the almost complex structure or the Kähler form) is an affine subspace of $\mathcal{SO}$ with vector space $\Gamma(T^*M \otimes \mathfrak u(M))$, where $\mathfrak u (M)$ is the subbundle of $\mathfrak{so}(M)$ formed by the endomorphisms which commute with $J$ (or in other words, the adjoint bundle of the $U(n)$-structure). Finally, we denote by $\mathfrak u(M)^{\perp}$ the orthogonal complement of $\mathfrak u(M)$ in $\mathfrak{so}(M)$, identified with the bundle of skew-symmetric endomorphisms of $TM$, anti-commuting with $J$.
\begin{defi}
The canonical Hermitian connection $\nb$ is the projection of $\nabla \in \mathcal{SO}$ on $\mathcal U$. Equivalently, it is the unique Hermitian connection such that $\nabla - \nb$ is a 1-form with values in $\mathfrak u(M)^{\perp}$.
\end{defi}

The difference $\eta = \nabla - \nb$ is known explicitely:
\[ \forall X \in TM, \quad \eta_X = \frac{1}{2} J \circ (\nabla_X J) \]
It measures the failure of the $U(n)$-structure to admit a torsion-free connection, in other words its \emph{torsion} or its 1-jet (see \cite{br}). Thus, it can be used (or $\nabla \omega$, or $\nabla J$) to classify almost Hermitian manifolds as in \cite{gr4}. For example, Kähler manifolds are defined by $\nabla$ itself being a Hermitian connection: $\nabla = \nb$. Equivalently, because $\eta$ determines $d\omega$ and the Nijenhuis tensor $N$, $\omega$ is closed and $J$ is integrable.

\begin{defi}
Let $M$ be an almost Hermitian manifold. The following conditions are equivalent and define a nearly Kähler manifold: \\
(\romannumeral 1) the torsion of $\nb$ is totally skew-symmetric.\\
(\romannumeral 2) $\forall X \in TM$, $(\nabla_X J)X = 0$ \\
(\romannumeral 3) $\forall X \in TM$, $\nabla_X \omega = \frac{1}{3} \iota_X d\omega$ \\
(\romannumeral 4) $d\omega$ is of type (3,0)+(0,3) and $N$ is totally skew-symmetric. \\
\label{nearly Kahler}
\end{defi}

The following result, due to Kirichenko \cite{ki}, is the base of the partial classification by Nagy in \cite{na2} of nearly Kähler manifolds.
\begin{prop}[Kirichenko]
For a nearly Kähler manifold, the torsion of the intrinsic connection is totally skew-symmetric (by definition) and parallel:
\[ \nb \eta = 0 \]
\label{parallel torsion}
\end{prop}
Moreover, this is equivalent to $\nb \nabla \omega = 0$ or $\nb d\omega = 0$ because $\nb$ is Hermitian.

\vs

Now, suppose that the curvature of $\nb$ is also parallel: $\nb \, \Rb = 0$. Then, $M$ is \emph{locally homogeneous} or an \emph{Ambrose-Singer manifold}. Besides, the associated \emph{infinitesimal model} is always \emph{regular} (for a definition of these notions, see \cite{cl2}, citing \cite{tr2}) and thus, if it is simply connected, $M$ is an homogeneous space. Examples obtained in this way belong to a particular class of homogeneous manifolds: the \emph{3-symmetric spaces}, defined by Gray \cite{gr2}, which shall interest us in the rest of this section. As expected, the 3-symmetric spaces are a generalization of the well-known \emph{symmetric} spaces:

\begin{defi}
A 3-symmetric space is a homogeneous space $M=G/H$, where $G$ has an automorphism $s$ of order 3 (instead of an involution, for a symmetric space) such that
\begin{equation} G^s_0 \subset H \subset G^s \label{H fixed points} \end{equation}
where $G^s=\{ g \in G \ | \ s(g)=g \}$ is the fixed points set of $s$ and $G^s_0$ is the identity component of $G^s$.
\end{defi}

Let $\g$, $\h$ be the Lie algebras of $G$, $H$, respectively. For a symmetric space, the eigenspace $\m$ for the eigenvalue $-1$ of $ds: \g \to \g$ (the derivative of $s$) is an ${\rm Ad}(H)$-invariant complement of $\h$ in $\g$, so that symmetric spaces are always {\it reductive}. Conversely, for a reductive homogeneous space, we define an endormorphism $f$ of $\g$ by $f|_{\h} = Id|_{\h}$ and $f|_{\m} = -Id|_{\m}$, which integrates into a Lie group automorphism of $G$ if and only if 
\begin{equation} 
[\m,\m] \subset \h \label{symmetric}
\end{equation}

For a 3-symmetric space, things are slightly more complicated because $ds$ has now three eigenvalues, $1$, $j=-\frac{1}{2} + i\frac{\sqrt{3}}{2}$ and $j^2=\bar j = -\frac{1}{2} - i\frac{\sqrt{3}}{2}$, two of which are complex. The corresponding eigenspaces decomposition is
\[ g^{\CM} = \h^{\CM} \oplus \m_j \oplus \m_{j^2}. \]
Setting $\m = (\m_j \oplus \m_{j^2}) \cap \g$, we get
\begin{equation} \g = \h \oplus \m, \quad {\rm Ad}(H)\m \subset \m \label{reductive} \end{equation}
so 3-symmetric spaces are also reductive. As a consequence, invariant tensors, for the left action of $G$ on $M$, are represented by constant, ${\rm Ad}(H)$-invariant tensors on $\m$. For example, an invariant almost complex structure is identified, first with the subbundle $T^+M \subset T^{\CM} M$ of (1,0)-vectors, then with a decomposition 

\begin{equation}
\m^{\CM} = \m^+ \oplus \m^- \quad \text{where} \quad {\rm Ad}(H)\m^+ \subset \m^+ \text{ and } \m^- = \overline{\m^+}.
\label{m+ + m-}
\end{equation}

\begin{defi}
The canonical almost complex structure of a 3-symmetric space is the invariant almost complex structure associated to $\m_j \oplus \m_{j^2}$.
\end{defi}

In other words, the restriction $ds: \m \to \m$ represents an invariant tensor $S$ of $M$ satisfying: \\
$(\romannumeral 1)$ $S^3 = Id$. \\
$(\romannumeral 2)$ For all $x \in M$, 1 is not an eigenvalue of $S_x$. \\
One can then write $S$ as for a (non trivial) third root of unity:
\begin{equation} 
S = -\frac{1}{2}Id + \frac{\sqrt{3}}{2}J
\label{J cubic root}
\end{equation}
where $J$ is the canonical almost complex structure of $M$.

Similarly, an ${\rm Ad}(H)$-invariant scalar product $g$ on $\m$ defines an invariant metric on $M$, also denoted by $g$, and the pair $(M,g)$ is called a {\it Riemannian} 3-symmetric space if and only $g$, $J$ are compatible.

Conversely a decomposition like (\ref{m+ + m-}) comes from an automorphism of order 3 if and only if $\h$, $\m^+$, $\m^-$ satisfy
\begin{equation} 
[ \m^+,\m^+ ] \subset \m^-, \quad [ \m^-,\m^- ] \subset \m^+ \quad \text{and} \quad [ \m^+,\m^- ] \subset \h^{\CM}.
\label{3symmetric}
\end{equation}
instead of (\ref{symmetric}). Now, conditions involving the Lie bracket might be interpreted, on a reductive homogeneous space, as conditions on the torsion and the curvature of the normal connection $\widehat \nabla$. The latter is defined as the $H$-connection on $G$ whose horizontal distribution is $G \times \m \subset TG \simeq G \times \g$.

\begin{lemm}
The torsion $\widehat T$ and the curvature $\widehat R$ of the normal connection $\widehat \nabla$, viewed as constant tensors, are respectively the $\m$-valued 2-form and the $\h$-valued 2-form on $\m$ given by
\[ \forall X,Y \in \m, \quad \widehat T(X,Y)=-[X,Y]^{\m}, \quad \widehat R_{X,Y} = [X,Y]^\h \]
\label{normal torsion curvature}
\end{lemm}

\begin{prop}
A reductive almost Hermitian homogeneous space $M=G/H$ is a 3-symmetric space if and only if it is quasi-Kähler and the intrinsic connection $\nb$ coïncides with $\widehat \nabla$.
\label{connexion normale = intrinseque}
\end{prop}

An almost Hermitian manifold is said \emph{quasi-Kähler} or \emph{(2,1)-symplectic} iff $d\omega$ has type (3,0)+(0,3) or equivalently $\eta$ (or $\nabla J$) is a section of $\LLambda^1 \otimes \mathfrak{u}(M)^{\perp} \cap \leftr \llambda^{2,0} \rightr \otimes TM$ where $\leftr \llambda^{2,0} \rightr \subset \LLambda^2$ is the bundle of real 2-forms of type (2,0)+(0,2).

\begin{proof}
By lemma \ref{normal torsion curvature}, (\ref{3symmetric}) is equivalent to
\[
\widehat T(\m^+,\m^+) \subset \m^-, \quad \widehat T(\m^-,\m^-) \subset \m^+, \quad \widehat T(\m^+,\m^-) = \{0\} 
\]
\[
\widehat R(\m^+,\m^+) = \widehat R(\m^-,\m^-) = \{ 0 \}
\]
The first line implies that $\widehat \eta = \nabla - \widehat \nabla \in \LLambda^1 \otimes \mathfrak{u}(M)^{\perp}$ (for a metric connection, the torsion and the difference with the Levi-Civita connection are in one-to-one correspondance), i.e. $\widehat \nabla$ is the canonical Hermitian connection. It also implies that $\eta = \widehat \eta \in \leftr \llambda^{2,0} \rightr \otimes TM$ (note that the bundles $\mathfrak{u}(M)^{\perp} \simeq \leftr \llambda^{2,0} \rightr$ are isomorphic through the operation of raising, or lowering, indices).
The second line is automatically satisfied for a quasi-Kähler manifold, see \cite{fa}.
\end{proof}

There is also a local version of proposition \ref{connexion normale = intrinseque}, announced in the middle of this section.

\begin{theo'}
An almost Hermitian manifold $M$ is locally 3-symmetric if and only if it is quasi-Kähler and the torsion and the curvature of the intrinsic connection $\nb$ satisfy
\[ 
\nb T = 0 \quad \text{and} \quad \nb \, \Rb = 0 
\]
\label{loc 3symmetric}
\end{theo'}

The definition of a {\em locally 3-symmetric space} given in \cite{gr2} relates to the existence of a family of \emph{local cubic isometries} $(s_x)_{x \in M}$ such that, for all $x \in M$, $x$ is an isolated fixed point of $s_x$ (for a 3-symmetric space, the automorphism $s$ provides such a family, moreover the isometries are globally defined).

Then, the requirement that $M$, the Riemannian 3-symmetric space, is a nearly Kähler manifold (which is more restrictive than quasi-Kähler) translates to a structural condition on the homogeneous space.

\begin{defi}
A reductive Riemannian homogeneous space is called naturally reductive iff the scalar product $g$ on $\m$ representing the metric satisfies
\[ 
\forall X,Y,Z \in \m, \quad g([X,Y],Z)=-g([X,Z],Y) 
\]
\end{defi}
Equivalently, the torsion $\widehat T$ of the normal connection is totally skew-symmetric. Now, for a 3-symmetric space, the intrinsic connection coincides with the normal connection, by proposition \ref{connexion normale = intrinseque}.
\begin{prop}
A Riemannian 3-symmetric space equipped with its canonical almost complex structure is nearly Kähler if and only if it is naturally reductive.
\label{3symmetric NK}
\end{prop}

\begin{NB}
Let $M=G/H$ be a compact, inner 3-symmetric space such that $G$ is compact, simple. The Killing form $B$ of $G$ is negative definite so it induces a scalar product $q=-B$ on $\g$. Then, the summand $\m$, associated to the eigendecomposition of $ds = {\rm Ad}(h)$, is orthogonal to $\h$ and the restriction of $q$ to $\m$ defines a naturally reductive metric that makes $M$ a nearly Kähler manifold.
\end{NB}

\section[$S^3 \times S^3$]{The case of $S^3 \times S^3$: the natural reduction to $SU(3)$}

In \cite{le}, Ledger and Obata gave a procedure to construct a nearly Kähler 3-symmetric space for each compact Lie group $G$. The Riemannian product $G \times G \times G$ has an obvious automorphism of order 3, given by the cyclic permutation, whose fixed point set is the diagonal subgroup $\Delta G = \{(x,x,x) \ | \ x \in G \} \simeq G$. The resulting homogeneous space is naturally isomorphic to $G \times G$: to fix things, we shall identify $(x,y)$ with $[x,y,1]$. In other words, we get a 3-symmetric structure on $G \times G$, invariant for the action:
\[
((h_1,h_2,h_3),(x,y)) \mapsto (h_1xh_3^{-1},h_2yh_3^{-1})
\]
Now, let $q$ be an ${\rm Ad}(G)$-invariant scalar product on $\g$, representing a biinvariant metric on $G$. We choose, for the ${\rm Ad}(\Delta G)$-invariant complement of $\delta \g$ (the Lie algebra of $\Delta G$) in $\g \oplus \g \oplus \g$, $\m = \{0\} \oplus \g \oplus \g$, the sum of the \emph{last} two factors, so that the restriction of $q \oplus q \oplus q$ to $\m$ defines a naturally reductive metric $g$ which is \emph{not} the biinvariant metric of $G \times G$. Indeed, the vector $(X,Y) \in \g \oplus \g \simeq T_e(G \times G)$ is identified with $(0,Y-X,X) \in \m$ so we have the explicit formula $g_e((X,Y),(X',Y')) = q(Y-X,Y'-X') + q(X,X')$.

Following this procedure and setting $G = SU(2) \simeq S^3$, we obtain a 6-dimensional example: $S^3 \times S^3$. In this section we will look for nearly Kähler structures on $S^3 \times S^3$ invariant for the smaller group 
\beq 
SU(2) \times SU(2) & \hookrightarrow & SU(2) \times SU(2) \times SU(2) \\
(h_1,h_2) & \mapsto & (h_1,h_2,1)
\eeq
or for the left action of $SU(2) \times SU(2)$ on itself. Such a structure is then simply given by constant tensors on the Lie algebra $\mathfrak{su}(2) \oplus \mathfrak{su}(2)$. The latter is a 6-dimensional vector space so, for example, the candidates for the metric belong to a 21-dimensional space $S^2(2\mathfrak{su}(2))$. Thus, the calculations involving the Levi-Civita connection as in definition \ref{nearly Kahler} are too hard and we shall look for another strategy.

\vs

Nearly Kähler manifolds in dimension 6 are special. In particular they're always Einstein \cite{gr3} and possess a Killing spinor \cite{gru}. But the most important feature, for us, is a natural reduction to $SU(3)$. Indeed, an $SU(3)$-structure in dimension 6, unlike an almost Hermitian structure, may be defined, without the metric, only by means of differential forms. This should be compared to the fact that a $G_2$-structure is determined by a differential 3-form on a 7-manifold.

In first approach, an $SU(3)$-manifold is an almost Hermitian manifold with a complex volume form (a complex 3-form of type (3,0) and constant norm) $\Psi$. This complex 3-form can be decomposed into real and imaginary parts: $\Psi = \psi + i\phi$ where
\begin{equation}
\forall X \in TM, \quad \iota_X \psi = \iota_{JX} \phi 
\label{phi=Jpsi} 
\end{equation}
As a consequence, one of the \emph{real} 3-forms $\psi$ or $\phi$, together with $J$, determines the reduction of the manifold to $SU(3)$. Now, for a nearly Kähler manifold, such a form is naturally given by the differential of the Kähler form. Indeed, because of our preliminaries, $d\omega$ has type (3,0)+(0,3) (see definition \ref{nearly Kahler}) and constant norm, since it is parallel for the metric connection $\nb$. 
\begin{defi}
The natural $SU(3)$-structure of a 6-dimensional nearly Kähler manifold is defined by
\[ \psi := \frac{1}{3} d\omega \]
where $\omega$ is the Kähler form.
\label{SU(3)-structure}
\end{defi}

There is more. Hitchin has observed in \cite{hi3} that a differential 2-form $\omega$ and a differential 3-form $\psi$, satisfying algebraic properties, are enough to define a reduction of the manifold to $SU(3)$. In particular they determine the metric $g$ and the almost complex structure $J$. Indeed, $SU(3)$ may be seen as the intersection of two groups, $Sp(3,\RM)$ and $SL(3,\CM)$, which are themselves the stabilizers of two exteriors forms on $\RM^6$. For the symplectic group $Sp(3,\RM)$, it is a non-degenerate 2-form of course. As for the second group, $GL(6,\RM)$ has two open orbits $\mathcal O_1$ and $\mathcal O_2$ on $\Lambda^3 \RM^6$. The stabilizer of the forms in the first orbit is $SL(3,\CM)$. To define the action of the latter, we must see $\RM^6$ as the complex vector space $\CM^3$. Consequently, if a differential 3-form $\psi$ belongs to $\mathcal O_1$ at each point, it determines an almost complex structure $J$ on $M$. Then, $\omega$, $J$ determine $g$ under certain conditions. 

We need to write this explicitly (after \cite{hi,hi3}). For $\psi \in \LLambda^3$, we define $K \in {\rm End}(TM) \otimes \LLambda^6$ by
\[
K(X) = A(\iota_X \psi \wedge \psi)
\]
where $A: \LLambda^5 \to TM \otimes \LLambda^6$ is the isomorphism induced by the exterior product. Then, $\tau(\psi)= \frac{1}{6} \mathrm{tr}\, K^2$ is a section of $(\LLambda^6)^2$ and it can be shown that
\[ K^2 = {\rm Id} \otimes \tau(\psi) \]
The 3-form $\psi$ belongs to $\mathcal O_1$ at each point if and only if
\begin{equation} \tau(\psi) < 0 \label{psi stable} \end{equation}
Then, for $\kappa = \sqrt{-\tau(\psi)}$,
\[ J = \frac{1}{\kappa} K \]
is an almost complex structure on $M$. Moreover, the 2-form $\omega$ is of type $(1,1)$ with respect to $J$ if and only if
\begin{equation} \omega \wedge \psi = 0 \label{omega (1,1)} \end{equation}
Finally, $\omega$ has to be non degenerate
\begin{equation} \omega \wedge \omega \wedge \omega \neq 0 \label{omega stable} \end{equation}
and $g$ has to be positive:
\begin{equation}
(X,Y) \mapsto g(X,Y) =
\omega(X,JY) > 0
\label{g positive}
\end{equation}

This is at the algebraic level. At the geometric level, Salamon and Chiossi \cite{ch} have shown that the 1-jet of the $SU(3)$-structure (or the intrinsic torsion) is completely determined by the differentials of $\omega$, $\psi$, $\phi$. In particular, nearly Kähler manifolds are viewed, in this section, as $SU(3)$-manifolds satisfying a first order condition. Thus, they're characterized by a differential system involving these three forms: 
\begin{equation}
\left\{
\begin{array}{ccc}
\psi & = & 3d\omega \\
d\phi & = & -2\mu \, \omega \wedge \omega
\end{array}
\right.
\label{diff system}
\end{equation}
where $\mu \in \RM$. This differential system was first written by Reyes Carri\'on in \cite{re}.
%The first line is the definition of the natural reduction to $SU(3)$ and implies $d\psi = 0$. Now $0=(d\Psi)^{1,3}=(d\psi)^{1,3} + i (d\phi)^{1,3}$ so the (3,1)+(1,3)-part of $d\phi$ also vanishes. The (2,2)-part may be written
%\[ (d\phi)^{2,2} = \mu \, \omega \wedge \omega + \nu \wedge \omega \]
%where $\mu$ is determined by the skew-symmetric part $N_a$ of $N$ and orthogonal part and $\nu \in [\llambda^{1,1}]$ is determined by the orthogonal part $N^{\perp}$ of $N_a$.

As a consequence, looking for a nearly Kähler structure on a manifold is the same as looking for a pair of forms $(\omega,\psi)$ satisfying (\ref{psi stable})-(\ref{g positive}) together with (\ref{diff system}) or, considering the particular form of (\ref{diff system}), for a 2-form $\omega$ only, satisfying a highly non linear second order differential equation. We shall resolve this system on the space of invariant 2-forms of $S^3 \times S^3$.

\vs

We work with a class of co-frames $(e_1,e_2,e_3,f_1,f_2,f_3)$, called \emph{cyclic co-frames}, satisfying: \\
(\romannumeral 1) $(e_1,e_2,e_3,f_1,f_2,f_3)$ is invariant for the action of $SU(2) \times SU(2)$ on itself.
(\romannumeral 2) the 1-forms $e_i$ (resp. $f_i$), $i=1,2,3$ vanish on the tangent space of the first (resp. the second) factor. \\
(\romannumeral 3) $de_i = e_{i+1} \wedge e_{i+2}$ where the subscripts are viewed as elements of $\ZM_3$. Similarly, $df_i = f_{i+1} \wedge f_{i+2}$.\\

The group of isometries of the sphere $S^3 \simeq SU(2)$, equipped with its round biinvariant metric is $SO(4)$, with isotropy subgroup $SO(3)$. Moreover the isotropy representation lifts to the adjoint representation of $SU(2) \simeq Spin(3)$. We denote $(u,l) \mapsto u.l$ the action of $SO(3)$ on the \emph{dual} $\mathfrak{su}(2)^*$ of the Lie algebra. Then, $SO(3) \times SO(3)$ acts transitively on the set of cyclic co-frames by
\begin{equation}
(u,v).(e_1,e_2,e_3,f_1,f_2,f_3) \mapsto (u.e_1,u.e_2,u.e_3,v.f_1,v.f_2,v.f_3) 
\label{change of frame}
\end{equation}
In other words, two such co-frames are exchanged by a diffeomorphism of $S^3 \times S^3$ and more precisely by an isometry of the canonical metric.

Now, a generic invariant 2-form may be written in the form
\begin{equation}
\omega = \sum_{i=1}^3 a_i e_{i+1} \wedge e_{i+2} + \sum_{i=1}^3 b_i
f_{i+1} \wedge f_{i+2} + \sum_{i,j=1}^3 c_{i,j}
e_i \wedge f_j
\label{omega=A+B+C}
\end{equation}
Let $A$ be the column vector of the $a_i$, $B$ the column vector of the $b_i$ and $C$ the square matrix $(c_{i,j})_{i,j = 1,2,3}$. The latter is subject to the following transformation rule in a change of cyclic co-frame (\ref{change of frame}): 
\begin{equation} C \mapsto MC ^tN, \label{change of frame A,B,C} \end{equation}
where the $3 \times 3$ matrices $M$ (resp. $N$) represent $u$ (resp. $v$) in the old base.

We have the first essential simplification of (\ref{omega=A+B+C}):

\begin{lemm}
Let $\omega$ be a non degenerate invariant 2-form on $S^3 \times S^3$. Then, $\omega$, $\psi = \frac{1}{3}d\omega$ satisfy (\ref{omega (1,1)}) if and only if there exists a cyclic co-frame $(e_1,e_2,e_3,f_1,f_2,f_3)$ such that
\begin{equation}
\omega = \lambda_1 \, e_1 \wedge f_1 + \lambda_2 \, e_2 \wedge f_2 + 
\lambda_3 \, e_3 \wedge f_3
\label{omega=lambda1,2,3}
\end{equation}
where $\forall i=1,2,3$, $\lambda_i \in \RM^*$.
\end{lemm}

\begin{proof}
Starting from (\ref{omega=A+B+C}), we calculate $\omega \wedge \omega \wedge \omega$. The 2-form $\omega$ is non degenerate if and only if
\begin{equation}
^tACB + {\rm det} \, C \neq 0 
\label{omega stable A,B,C}
\end{equation}
Then, we calculate $\psi = \frac{1}{3}d\omega$ using the relations (\romannumeral 3), in the definition of a cyclic co-frame. We find that $\omega \wedge \psi = 0$ is equivalent to
$ ^tAC = CB = 0$. Reintroducing these equations in (\ref{omega stable A,B,C}), we get ${\rm det}\, C = 0$, i.e. $C$ is nonsingular, and so $A = B = 0$. 

Secondly, we can always suppose that $C$ is diagonal. Indeed, we write $C$ as the product of a symmetric matrix $S$ and an orthogonal matrix $O \in SO(3)$. We diagonalize $S$: there exists an orthogonal matrix $P$ such that $S = \, ^t\! PDP$ where $D$ is diagonal. Thus $C = \, ^t\!P D (PO)$ and by (\ref{change of frame A,B,C}), $\omega$ can always be written in the form (\ref{omega=lambda1,2,3}) where $D = {\rm diag}(\lambda_1,\lambda_2,\lambda_3)$.
\end{proof}

This is a key lemma that will constitute the base of our next calculations.
\begin{lemm}
Let $\omega$ be the invariant 2-form given by (\ref{omega=lambda1,2,3}) in a cyclic co-frame. Then, $\omega$, $\psi = \frac{1}{3} d\omega$ define an $SU(3)$ structure on $S^3 \times S^3$ if and only if: \\
(\romannumeral 1) $(\lambda_1 - \lambda_2 - \lambda_3)(-\lambda_1 + \lambda_2 - \lambda_3)(-\lambda_1 - \lambda_2 + \lambda_3)(\lambda_1 + \lambda_2 + \lambda_3) < 0$ \\
(\romannumeral 2) $\lambda_1 \lambda_2 \lambda_3 > 0$
\label{omega,psi SU(3) lambda1,2,3}
\end{lemm}

\begin{proof}
The first condition is simply (\ref{psi stable}). Indeed,
\[ 3\psi = \lambda_1(e_{23} \wedge f_1 - e_1 \wedge f_{23}) + \lambda_2(e_{31} \wedge f_2 - e_2 \wedge f_{31}) + \lambda_3 (e_{12} \wedge f_3 - e_3 \wedge f_{12}) \]
and we calculate
\begin{equation}
81\tau(\psi) = (\lambda_1^4 + \lambda_2^4 + \lambda_3^4
- 2\lambda_1^2\lambda_2^2 - 2\lambda_2^2\lambda_3^2 - 2\lambda_1^2\lambda_3^2) \otimes {\rm vol}^2 
\label{psi in O1 ?}
\end{equation}
where ${\rm vol} = e_{123} \wedge f_{123}$. Now, the polynomial of degree 4 on the $\lambda_i$ in (\ref{psi in O1 ?}) factors into (\romannumeral 1) of lemma \ref{omega,psi SU(3) lambda1,2,3}.

The second condition comes from the positivity of the metric. Note that the product $\lambda_1 \lambda_2 \lambda_3 = {\rm det}\, C$ is independent on the choice of the co-frame $(e_1,e_2,e_3,f_1,f_2,f_3)$, the determinant of the matrices $M$, $N$, in (\ref{change of frame A,B,C}), being equal to 1. First, we compute the almost complex structure in the dual frame $(X_1,X_2,X_3,Y_1,Y_2,Y_3)$:
\begin{equation}
JX_i= \alpha_i X_i + \beta_i Y_i, \quad JY_i = -\beta_i X_i - \alpha_i Y_i
\label{J lambda1,2,3}
\end{equation}
where
\[ \alpha_i = \frac{1}{k}(\lambda_i^2 - \lambda_{i+1}^2 - \lambda_{i+2}^2) \quad \text{and} \quad \quad \beta_i = -\frac{2}{k}\lambda_{i+1}\lambda_{i+2} \]
\[ 9\kappa = k \, {\rm vol} \]
Now, $X \mapsto \omega(X,JX)$ is the sum of three quadratic forms of degree 2: 
\[ q_i: (x_i,y_i) \mapsto \frac{2\lambda_i}{k}\big(\lambda_{i+1}\lambda_{i+2}x_i^2-
(\lambda_i^2-\lambda_{i+1}^2-\lambda_{i+2}^2)x_i y_i + 
\lambda_{i+1}\lambda_{i+2}y_i^2\big)
\]
The discriminant of these forms is $k^2 > 0$ so they are definite and their sign is given by $\lambda_1 \lambda_2 \lambda_3$.
\end{proof}

Using (\ref{J lambda1,2,3}) it is easy to compute $\phi$, by (\ref{phi=Jpsi}), and translate the second line of (\ref{diff system}) into a system of equations \emph{on the $\lambda_i$}. We refer to \cite{bu} for a detailed proof.

\begin{lemm}
Let $\omega$ be the invariant 2-form given by (\ref{omega=lambda1,2,3}) in a cyclic co-frame. Then $\omega$, $\psi = \frac{1}{3}d\omega$ induce a nearly Kähler structure on $S^3 \times S^3$ if and only if there exists $\mu \in \RM$ such that, for all $i=1,2,3$,
\begin{equation}
c =\lambda_i^2(\lambda_i^2-\lambda_{i+1}^2-\lambda_{i+2}^2)
\label{diff system lambda1,2,3}
\end{equation}
where 
\[ c = -2\mu k \, {\rm det}\,C. \]
\end{lemm}

Finally we can conclude:

\begin{prop}
There exists a unique (up to homothety, up to a sign) left-invariant nearly Kähler structure on $S^3 \times S^3$, corresponding to Ledger and Obata's construction of a 3-symmetric space.
\label{S3xS3}
\end{prop}

\begin{proof}
Thanks to the preparatory work, we only need to solve the system (\ref{diff system lambda1,2,3}) of equations of degree 4 on the $\lambda_i$.  Let $\Lambda = \lambda^2_1 + \lambda^2_2 + \lambda^2_3$. For all $i = 1,2,3$, $\lambda_i^2$ is a solution of the unique equation of degree 2:
\begin{equation}
2x^2 - \Lambda x - c = 0 
\label{equation °2}
\end{equation}
Suppose that $\lambda^2_1$, $\lambda^2_2$ are two distinct solutions of (\ref{equation °2}) and $\lambda^2_3=\lambda^2_2$, for example. The sum of the roots $\lambda^2_1 + \lambda^2_2$ equals $\displaystyle{\frac{\Lambda}{2}}$. But then, we also have $\Lambda = \lambda^2_1 + 2\lambda^2_2$, by definition of $\Lambda$. We immediatly get $\lambda_1=0$, i.e. $\omega$ is degenerate, a contradiction. Thus, the $\lambda_i$ must be equal, up to a sign. The positivity of the metric, (\ref{g positive}) or \ref{omega,psi SU(3) lambda1,2,3} (\romannumeral 2), implies that the three signs are positive, or only one of them. These two solutions are in fact the same: one is obtained from the other by a rotation of angle $\pi$ in the first factor. Finally, one can always write $\omega$, for a left-invariant nearly Kähler structure, 
\[ \omega = \lambda(e_1 \wedge f_1 + e_2 \wedge f_2 + e_3 \wedge f_3) \] 
where $\lambda = \lambda_1 = \lambda_2 = \lambda_3 \in \RM^{+*}$. 

We also have $k = \lambda^2\sqrt{3}$,
\begin{equation}
JX_i = \frac{1}{\lambda^2\sqrt{3}} ( -X_i + 2Y_i), \quad  JY_i = \frac{1}{\lambda^2\sqrt{3}}( -2X_i + Y_i). 
\label{J S3xS3}
\end{equation}
Now, (\ref{J S3xS3}) coincides with the canonical almost complex structure of $SU(2) \times SU(2) \times SU(2)/\Delta SU(2)$. Indeed, the automorphism of order 3, $s : (h_1,h_2,h_3) \mapsto (h_2,h_3,h_1)$, induces the endomorphism $S : (X,Y) \mapsto (Y-X,-X)$ on $\m \simeq \mathfrak{su}(2) \oplus \mathfrak{su}(2)$. Then, by (\ref{J cubic root}), $J$ is identified with
\[ J : (X,Y) \mapsto \frac{1}{\sqrt{3}}(2Y-X,-2X+Y) \]
This is no more than (\ref{J S3xS3}) with $\lambda = 1$ for an appropriate choice of base. \\
NB : $c=2\mu \lambda^5 \sqrt{3}$ so by (\ref{diff system lambda1,2,3}), $\displaystyle{\mu = \frac{1}{2\lambda \sqrt{3}}}$ is inversely proportional to the norm of $\omega$.
\end{proof}

\section[Twistors spaces]{Twistors spaces: the complex projective space $\CM P^3$ and the flag manifold $\mathbb F^3$}

The twistor space $Z$ of a 4-dimensional, Riemannian, oriented manifold $(N,h)$ is equipped with two natural almost complex structures. %To describe them: the fibre $\CM P^1$ has a natural complex structure $J_0$ which may be completed, using the horizontal distribution $H$ of the Levi-Civita connection of the base, in two different manners, depending on a choice of orientation of the fibre. Let $\pi: Z \to N$ be the twistor fibration. A point $j \in Z$ is identified with an complex structure on $T_x N$, $x=\pi(j)$, compatible with the metric and orientation. Then, we define $J_+$, $J_-$ by
%\[ J_{\pm}|_{V_j} = \pm J_0, \quad J_{\pm}|_{H_j} = \pi^* j \]
%where $V \subset TZ$ is the vertical distribution associated with $\pi$ (NB: the differential $\pi^*$ of the projection induces an isomorphism between $H_j$ and $T_x N$). 
The first, $J_+$, studied by Atiyah, Hitchin and Singer \cite{at}, is integrable as soon as $N$ is \emph{self-dual}, i.e one half of the Weyl tensor of $h$ vanishes, while the second $J_-$, which was first considered by Salamon and Eels in \cite{ee}, is never integrable. Now, on $Z$, varying the scalar curvature of the fibre, there are also a 1-parameter family of metrics $(g_t)_{t \in ]0,+\infty[}$ such that the twistor fibration over $(N,h)$ is a Riemannian submersion and for all $t$, $(g_t,J_{\pm})$ is an almost Hermitian structure on $Z$. A natural problem is then to look at the type of this almost Hermitian structure. We are particularly interested in the case where $N$ is compact, self-dual and Einstein.

\begin{theo'}[Hitchin, Eels \& Salamon]
Let $(N,h)$ be a compact, Riemannian, oriented, self-dual, Einstein 4-manifold, $Z$ its twistor space and $(g_t)_{t \in [0,+\infty[}$, the twistor metrics. There exists a choice of parameter such that the scalar curvature of the fibre of $\pi: Z \to N$ is proportionnal to $t$ and \\
(\romannumeral 1) $(Z,g_2,J_+)$ is Kähler \\
(\romannumeral 2) $(Z,g_1,J_-)$ is nearly Kähler.
\label{Kähler-NK}
\end{theo'}

This provides us for two compact, homogeneous, nearly Kähler structures in dimension 6 on the complex projective space $\CM P^3$ and the flag manifold $\mathbb F^3$, the twistor spaces of $S^4$ and $\CM P^2$, respectively. Moreover, we shall see that they correspond to a 3-symmetric structure. The goal of this section is to prove that these are the only invariant nearly Kähler structures on the above mentioned spaces.

\vs

This is quite easy for $\CM P^3$. The complex projective space is seen, in this context, as $Sp(2)/U(1)Sp(1)$. More generally, $\CM P^{2q+1}$ is isomorphic to $Sp(q+1)/U(1)Sp(q)$. Indeed, $Sp(q+1)$ acts transitively on $\CM^{2q+2} \simeq \HM^{q+1}$, preserving the complex lines, and the isotropy subgroup at $x \in \CM P^{2q+1}$ fixes also $jx$ and acts on the orthogonal of $\{x, jx \}$, identified with $\HM^q$, as $Sp(q)$. Representing $\mathfrak{sp}(q)$ as usual in the set of $q \times q$ matrices, the embedding of $\h = \mathfrak{u}(1) \oplus \mathfrak{sp}(1)$, the Lie algebra of $H=U(1)Sp(1)$, in $\g=\mathfrak{sp}(2)$ is given by the composition of the natural maps $\mathfrak{u}(1) \hookrightarrow \mathfrak{sp}(1)$ and
the identity of $\mathfrak{sp}(1)$, followed by the diagonal map $\mathfrak{sp}(1) \oplus  \mathfrak{sp}(1) \hookrightarrow
\mathfrak{sp}(2)$. Thus, a natural choice of complement of $\h$ in $\g$ is $\m = \mathfrak{p} \oplus \mathfrak{v}$ where
\[
\mathfrak p=\{
\left(\begin{array}{cc}
0   & a \\
a^* & 0
\end{array}\right)
| \ a \in \HM \}
\quad \text{ and } \quad
\mathfrak v = \{
\left(\begin{array}{cc}
b   & 0\\
0 & 0
\end{array}\right)
| \ b = jx + ky, \ x,y \in \RM \}
\]
These two subsets are ${\rm Ad}(H)$-invariant so their sum is too and $\m$ may be identified with the isotropy representation of $Sp(2)/U(1)Sp(1)$. The restriction of ${\rm Ad}(H)$ to $\mathfrak{p}$ is irreducible because the induced representation of $Sp(1)$ is isomorphic to the standard one on $\HM$. Similarly, the restriction of ${\rm Ad}(H)$ to $\mathfrak{v}$ induces the standard representation of $U(1)$ on $\CM$. As a consequence, $\m$ has exactly two irreducible summands and the set of invariant metrics on $\CM P^3$ has dimension 2. Moreover, we can gain one degree of freedom by working "up to homothety". Finally, we get a 1-parameter family of metrics which may be identified with $(g_t)_{t \in [0,+\infty[}$. On the other hand, $\CM P^3$ has $2^2=4$ invariant almost complex structures according to \cite{wo}, theorem 4.3: $\pm J_+$ and $\pm J_-$. Thus, we are in the hypothesis of Muskarov's work \cite{mu}. The conclusion is, as in proposition \ref{Kähler-NK},
\begin{prop}
The homogeneous space $\CM P^3 \simeq Sp(2)/U(1)Sp(1)$ has a unique -- up to homotethy, up to a sign -- invariant nearly Kähler structure $(g_1,J_-)$, associated to the twistor fibration of $S^4$.
\label{CP3}
\end{prop}

\vs

Things are more complicated for $\mathbb F^3$ because the isotropy representation has three irreducible summands. We could adapt Muskarov's proof \cite{mu} for this case, as was done in \cite{bu}, calculating $\nabla J$ for all $g$ with Levi-Civita connection $\nabla$ in the 3-dimensional space of invariants metrics, and all invariant almost complex structure $J$, or we can look more carefully at the structure of $\mathbb F^3$, as we will do now.

The flag manifold is the space of pairs $(l,p)$ where $p \subset \CM^3$ is a complex plane and $l \subset p$ is a complex line. It is isomorphic to $G/H$ where $G = U(3)$ and $H = S^1 \times S^1 \times S^1$. Indeed, we see $U(3)$ as the space of unitary bases of $\CM^3$. Then, the map $(u_1,u_2,u_3) \mapsto (l,p)$ where $l= \CM u_1$ and $p = \CM u_1 \oplus \CM u_2$, is an $H$-principal bundle over $\mathbb F^3$ with total space $G$. 

Now, the maps $\pi_a: (l,p) \mapsto \CM u_a$, $a=1,2,3$ are well-defined (because $\CM u_1 = l$ ; $\CM u_2 = l^{\perp}$ is the orthogonal line of $l$ in $p$ and finally $\CM u_3 = p^{\perp}$ is the orthogonal of $p$) and give three different $\CM P^1$-fibrations from $\mathbb F^3$ to $\CM P^2$ (or, if we identify the base spaces, three different realizations of $\mathbb F^3$ as an almost complex submanifold of $\mathcal Z$, the twistor space of $\CM P^2$). This has the following geometrical interpretation: on the fibre of $\pi_3$, $l$ varies inside $p$ ; on the fibre of $\pi_1$, it is the plane $p$ that varies around the line $l$ ; finally, on the fibre of $\pi_2$, both $l$ and $p$ vary while $l^{\perp}$ is fixed. These fibrations are all \emph{twistor} fibrations over an Einstein self-dual 4-manifold. Let $I_+$, $J_+$, $K_+$ and $I_-$, $J_-$, $K_-$ be the associated almost complex structures, as in proposition \ref{Kähler-NK}. What is remarkable is that the Kählerian structures are all distinct but the non-integrable nearly Kähler structures coïncide: $I_- = J_- = K_-$. This observation was already made by Salamon in \cite{sa2}, section 6. We shall prove this at the infinitesimal level. As a consequence, the four above almost complex structures, and their opposites, exhaust all the $8=2^3$ invariant almost complex structures on $\mathbb F^3 \simeq U(3)/(S^1)^3$. 

Let $\g=\mathfrak{u}(3)$ be the set of the trace-free, anti-Hermitian, $3 \times 3$ matrices. Then, $\h=3\mathfrak{u}(1)$ is identified with the subgroup of the diagonal matrices and the set $\m$ of the matrices with zeros on the diagonal is an obvious ${\rm Ad}(H)$-invariant complement of $\h$ in $\g$. Denote
\[
\la a,b,c \ra = 
\left(\begin{array}{ccc}
0 & -\overline{a} & b \\
a & 0 & -\overline{c} \\
-\overline{b} & c & 0
\end{array}\right)
\]
for all $a,b,c \in \CM$.
\begin{equation}
\forall h = 
\left(\begin{array}{ccc}
e^{ir} & 0 & 0 \\
0 & e^{is} & 0 \\
0 & 0 & e^{it}
\end{array}\right)
\in H, \ Ad_h \la a,b,c \ra = \la e^{i(s-t)}a,e^{i(t-r)}b,e^{i(r-s)}c \ra
\label{Ad(H)m flags}
\end{equation}
It is easily seen on (\ref{Ad(H)m flags}) that the isotropy representation decomposes into
\begin{equation}
\m = \mathfrak{p} \oplus \mathfrak{q} \oplus \mathfrak{r} 
\label{irr decomposition flags}
\end{equation}
where
\begin{eqnarray*}
\mathfrak{p} & = & \{\la a,0,0 \ra \mid a \in \CM \} \\
\mathfrak{q} & = & \{\la 0,b,0 \ra \mid b \in \CM \} \\
\mathfrak{r} & = & \{\la 0,0,c \ra \mid c \in \CM \}
\end{eqnarray*}

Each of these 2-dimensional subspaces, $\mathfrak{a}$, has a natural scalar product $g_{\mathfrak{a}}$ and a natural complex structure $J_{\mathfrak{a}}$. For example, on $\mathfrak p$,
\[ g_{\mathfrak{p}}(\la a,0,0 \ra,\la a',0,0 \ra) = Re(a\overline{a'}) \quad \text{and} \quad J_{\mathfrak{p}}\la a,0,0 \ra = \la ia,0,0 \ra \]
Moreover, we denote by $\mathfrak{p}^+ \oplus \mathfrak{p}^-$ the decomposition of $\mathfrak{p}_{\CM}$ associated to $J_{\mathfrak{p}}$ (and similarly for $\mathfrak{q}$, $\mathfrak{r}$). The relations between the three subspaces $\mathfrak{p}$, $\mathfrak{q}$ and $\mathfrak{r}$ are given by the Lie brackets:
\begin{equation}
[ \la a,0,0 \ra,\la 0,b,0 \ra ] = \la 0,0,ab \ra, \ \text{ etc.}
\label{brackets a,b,c flags}
\end{equation}
and
\begin{equation}
[\la a,0,0 \ra,\la a',0,0 \ra ] = {\rm diag}\,(iy,-iy,0) \in \h, \ \text{ where } y = 2{\rm Im}\,(a\overline{a'}),\  \text{ etc.}
\label{brackets a,a' flags}
\end{equation}
From (\ref{brackets a,b,c flags}) and (\ref{brackets a,a' flags}), it is easy to calculate
\[ [\mathfrak{p}^+,\mathfrak{q}^+] \subset \mathfrak{r}^-, \quad [\mathfrak{p}^+,\mathfrak{r}^-] \subset \mathfrak{q}^+ \quad \text{and} \quad [\mathfrak{p}^+,\mathfrak{p}^+] = \{0\}, \text{ etc.} \]
Thus, $\mathfrak{p}^+ \oplus \mathfrak{q}^+ \oplus \mathfrak{r}^-$ is a subalgebra of $\m_{\CM}$, corresponding to an invariant \emph{complex} structure on $\mathbb F^3$: $K_+$. In the same way, $I_+$, $J_+$ are represented by $\mathfrak{p}^- \oplus \mathfrak{q}^+ \oplus \mathfrak{r}^+$, $\mathfrak{p}^+ \oplus \mathfrak{q}^- \oplus \mathfrak{r}^+$, respectively. On the contrary, $\m^+ = \mathfrak{p}^+ \oplus \mathfrak{q}^+ \oplus \mathfrak{r}^+$ is \emph{not} a subalgebra. However, it satisfies (\ref{3symmetric}) with $\m^- = \overline{\m^+} = \mathfrak{p}^- \oplus \mathfrak{q}^- \oplus \mathfrak{r}^-$. Thus, $\mathbb F^3$ is a 3-symmetric space with $\m^+$ as canonical almost complex structure. By proposition \ref{connexion normale = intrinseque}, each pair $(g,\pm J_-)$, where $g$ is a generic invariant metric,
\[ g = rg_{\mathfrak{p}} + sg_{\mathfrak{q}} + tg_{\mathfrak{r}}, \]
is a (2,1)-symplectic homogeneous structure on $\mathbb F^3$. Moreover, since $I_+$, $J_+$, $K_+$ are integrable, every invariant (strictly) nearly Kähler structure on $\mathbb F^3$ has that form, where $g$ is naturally reductive (see proposition \ref{3symmetric NK}). It is not hard to see that this corresponds to $r=s=t$. 

\begin{prop}
The nearly Kähler structures associated to the three natural twistor fibrations over $\CM P^2$ on $\mathbb F^3 \simeq U(3)/(S^1)^3$ coincide. Moreover, every invariant nearly Kähler structure on $\mathbb F^3$ is proportional to this one (or to its opposite).
\label{F3}
\end{prop}

\begin{NB}
The decomposition (\ref{irr decomposition flags}) is still the irreducible decomposition for $SU(3)$ $\subset U(3)$, so the results remain valid for this smaller group of isometries. This observation will be useful in section 5.
\end{NB}

%The two structures coincide on the horizontal distribution $H \subset TZ$ associated to the Levi-Civita connection of $h$, and on the vertical distribution, parallel to the fibres, the almost complex structures have opposite sign and the Kähler metric is twice the nearly Kähler. This situation is more general. There exists a Kählerian but also a nearly Kähler structure similarly related first on the the twistor spaces of positive quaternion-Kähler manifolds. Thus, the latter can be considered as an analog of Einstein self-dual manifolds in dimension 4. This includes the twistor spaces of the Wolf spaces i.e. compact, symmetric, quaternion-Kähler manifolds. Now the construction can be adapted to a larger class of symmetric spaces. Moreover in this case the twistor space is 3-symmetric. 

\section[The 6-sphere]{Weak holonomy and special holonomy: the case of the sphere $S^6$}

While $S^3 \times S^3$ is considered the hardest case because the isotropy is reduced to $\{ 0 \}$, the case of $S^6$ is apparently the easiest because the isotropy is maximal, $H = SU(3)$. The isotropy representation is the standard representation of $SU(3)$ in dimension 6, in particular it is irreducible so there exists only one metric up to homothety and also one almost complex structure up to a sign preserved by $H$. However, a difficulty occurs since on $S^6$ unlike on the other spaces considered in sections 2 or 3, the metric doesn't determine the almost complex structure:
\begin{prop}
Let $(M,g)$ be a complete Riemannian manifold of dimension 6, not isomorphic to the round sphere. If there exists an  almost complex structure $J$ on $M$ such that $(M,g,J)$ is nearly Kähler (non Kählerian) then it is unique. Moreover, in this case, $J$ is invariant by the isometry group of $g$.
\end{prop}
This can be proved using the spinors (by \cite{gru}, a $6$-dimensional Riemannian manifold admits a nearly Kähler structure if and only if it carries a real Killing spinor): see \cite{bu} proposition 2.4 and the reference therein \cite{bau}, proposition 1, p126, or by a "cone argument" (see below) as in \cite{ve}, proposition 4.7.

On the contrary, on the sphere $S^6$ equipped with its round metric $g_0$, there exist infinitely many compatible nearly Kähler structures:
\begin{prop}
The set $\mathcal J$ of almost complex structures $J$ such that $(S^6,g_0,J)$ is nearly Kähler, is isomorphic to $SO(7)/G_2 \simeq \RM P^7$.
\label{NK structures on S6}
\end{prop}

\begin{coro}
Nearly Kähler structures compatible with the canonical metric on $S^6$ -- or almost complex structures $J \in \mathcal J$, are all conjugated by the isometry group $SO(7)$ of $g_0$.
\label{S6}
\end{coro}

To show this, we shall use a theorem of Bär \cite{ba}: the Riemannian cone of a nearly Kähler manifold has holonomy contained in $G_2$. However, in order to remain faithful to the point of view of differential forms adopted in this article we prefer the presentation by Hitchin \cite{hi3} of this fact. According to section 2, a nearly Kähler structure is determined by a pair of differential forms $(\omega,\psi)$ satisfying (\ref{diff system}) as well as algebraic conditions (\ref{psi stable})-(\ref{g positive}). Moreover there exists, around each point, an orthonormal co-frame $(e_1,\ldots,e_6)$ such that 
\[ \omega = e_{12} + e_{34} + e_{56}  \]
\[ \text{and} \quad \psi = e_{135} - e_{146} - e_{236} - e_{245} , \]
where $e_{12} = e_1 \wedge e_2$, $e_{135} = e_1 \wedge e_3 \wedge e_5$, etc. Now, the cone of $(M,g)$ is the Riemannian manifold $(\overline M, \overline g)$ where $\overline M = M \times \RM^+$ and $\overline g = r^2 g + dr^2$ in the coordinates $(x,r)$. We define a section $\rho$ of $\Lambda^3 \overline M$ by
\begin{equation}
\rho = r^2 dr \wedge \omega + r^3 \psi 
\label{3-form cone}
\end{equation}
Let $u_0 = dr$ and for all $i=1,\ldots,6$, $u_i=\frac{1}{r}e_i$.
\[ \rho =  u_{012} + u_{034} + u_{056} + u_{135} - u_{146} - u_{236} - u_{245}\]
Thus, $\rho$ is a generic 3-form, inducing a $G_2$-structure on the 7-manifold $\overline M$ such that $(u_0,\ldots,u_6)$ is an orthonormal co-frame for the underlying Riemannian structure i.e. $\overline g$ is the metric determined by $\rho$, given the inclusion of $G_2$ in $SO(7)$. Moreover, if we denote by $*$ the Hodge dual of $\overline g$,
\[ *\rho =  - r^3 dr \wedge \phi + \frac{\mu}{2} r^4 \omega \wedge \omega   \]
Then, (\ref{diff system}) implies
\begin{equation*}
\left\{
\begin{array}{ccc}
d\rho & = & 0 \\
d*\rho & = & 0
\end{array}
\right.
\end{equation*}
By Gray, Fernandez \cite{fe}, the last couple of equations is equivalent to $\nabla^{\overline g} \rho = 0$, where $\nabla^{\overline g}$ is the Levi-Civita connection of $\overline g$. In other words, the holonomy of $(\overline M,\overline g)$ is contained in $G_2$. 

Reciprocally, a parallel, generic 3-form on $\overline M$ can always be written (\ref{3-form cone}) where $(\omega,\psi)$ define a nearly Kähler $SU(3)$-structure on $M$.

We are now ready to prove proposition \ref{NK structures on S6}.
\begin{proof}
The Riemannian cone of the 6-sphere is the Euclidean space $\RM^7$. According to what precedes, nearly Kähler structure on $S^6$, compatible with $g_0$, define a parallel or equivalently, a constant 3-form on $\RM^7$. This form must have the appropriate algebraic type, i.e. be an element of the open orbit $\mathcal O \simeq GL(7,\RM)/G_2 \subset \Lambda^3 \RM^7$. But it must also induce the good metric (the cone metric) on $\RM^7$. Finally the 3-forms that parametrize $\mathcal J$ belong to a subset of $\mathcal O$, isomorphic to $SO(7)/G_2$.
\end{proof}

The homogeneous nearly Kähler structure on $S^6$ is defined using the octonions. The octonian product $(x,y) \mapsto x.y$ may be described in the following way. First, the 8-dimensional real vector space $\mathbb O$ decomposes into $\RM \oplus \Im$. The subspace $\Im \simeq \RM^7$ is called the space of imaginary octonions and equipped, for our purposes, with an inner product $(x,y) \mapsto \la x , y \ra$. Secondly, rules (\romannumeral 1) to (\romannumeral 4) below are satisfied with respect to this decomposition: \\
(\romannumeral 1) $1.1 = 1$, \\
(\romannumeral 2) for all $x \in \Im$, $1.x = x$, \\
(\romannumeral 3) for all imaginary quaternion $x$ of norm 1, $x.x = -1$, \\
(\romannumeral 4) finally, for all orthogonal $x,y \in \Im$, $x.y = P(x,y)$ where $P: \RM^7 \times \RM^7 \to \RM^7$ is the \emph{2-fold vector cross product}. The latter satisfies himself $(x,y,z) \mapsto \la P(x,y),z \ra$ is a 3-form. In particular $P(x,y)$ is orthogonal to $x$, $y$.

Now, let $S^6$ be the unit sphere in $\Im \simeq \RM^7$. The tangent space at $x \in S^6$ is identified with the subspace of $\RM^7$ orthogonal to $x$. Then, $J$ is defined by
\beq
J_x: T_x M & \to & T_x M \\
          y & \mapsto & x.y
\eeq
This is a well defined almost complex structure because $x.y$ is orthogonal to $x$ by rule number (\romannumeral 4) and $J^2=-Id$ by rules (\romannumeral 2), (\romannumeral 3). Moreover $J$ is compatible with the metric induced by $\la . , . \ra$ on $S^6$ (the round metric $g_0$). 

Reciprocally, starting from $J \in \mathcal J$, we may rebuild an octonian product on $\RM^7$ by bilinearity. Then $J$ is invariant for the associated group of automorphisms, isomorphic to $G_2$. Consequently, all nearly Kähler structures on the round sphere are of the same kind (homogeneous). They correspond to different choices of embeddings of $G_2$ into the group of isometries of $(S^6,g_0)$. These embeddings are parametrized by $SO(7)/G_2 \simeq \RM P(7)$ thus the above discussion gives an alternative proof of proposition \ref{NK structures on S6}.

The same question was raised by Friedrich in \cite{fr4}. His proof is very similar to our first one though it uses the Hodge Laplacian instead of the differential system (\ref{diff system}). Moreover, another proof is mentioned, that uses the Killing spinors.

Finally, a third or fifth way of looking at this is the following: each $SU(3)$-structure $(\omega_x,\psi_x)$ on a tangent space $T_x S^6$, $x \in S^6$, may be extended in a unique way to a nearly Kähler structure on the whole manifold. Let $\rho$ be the constant 3-form on $\RM^7$ whose value at $(x,1)$ is
\[ \rho_{(x,1)} = dr \wedge \omega_x + \psi_x \]
Then, $\rho$ is parallel for the Levi-Civita connection of the flat metric and in return, $\omega = \iota_{\partial_r} \rho$, $\psi=\frac{1}{3}d\omega$ determine a nearly Kähler structure on $S^6$ whose values at $x$ coincide with $\omega_x$, $\psi_x$, consistent with our notations.
Now, $SU(3)$-structures on a $6$-dimensional vector space are parametrized by $SO(6)/SU(3)$ so this constitutes a (differential) geometric proof of the isomorphism
\[ \frac{SO(6)}{SU(3)} \simeq \frac{SO(7)}{G_2} \simeq \RM P(7) \simeq \mathcal J \]

\section[Classification]{Classification of 3-symmetric spaces and proof of the theorems}

In this section, we draw all the useful conclusions of the facts gathered in the previous sections about nearly Kähler manifolds and synthesize all the results to achieve the proof of theorems \ref{classification} and \ref{resolution conjecture}.

\vs

As such, the conjecture of Gray and Wolf is not easy to settle because the notion of a nearly Kähler manifold, or even of a 3-symmetric space, are too rich. Indeed, the classification of 3-symmetric spaces \cite{wo} discriminates between three types, that correspond to quite different geometries:

\begin{itemize}

\item[$\mathfrak A$.] The first type consists in twistor spaces of symmetric spaces. Indeed, the situation described in the beginning of section 3 has a wider application than the 6-dimensional twistor spaces of Einstein self-dual 4-manifolds. First, the study of quaternion-Kähler manifolds, i.e. Riemannian manifolds whose holonomy is contained in $Sp(q)Sp(1)$, provides an analog of that situation in dimension $4q$, $q \geq 2$. Such manifolds still admit a twistor space $Z \to M$ with fibre $\CM P^1$ and two almost complex structures $J_+$ and $J_-$ related by
\[ J_+|_V = -J_-|_V, \quad  J_+|_H = J_-|_H, \]
where $V$ is the vertical distribution, tangent to the fibres, and $H$ is the horizontal distribution of the Levi-Civita connection of the base, such that the first one is always integrable (see for example \cite{sa3}) while the second, $J_-$, is non-integrable. Moreover, for a \emph{positive} quaternion-Kähler manifold, there exist two natural metrics $g_1$, $g_2$ such that $(g_2,J_+)$ is a Kählerian structure and $(g_1,J_-)$ is a nearly Kähler structure on $Z$ (compare with proposition \ref{Kähler-NK}). This includes the twistor spaces of the Wolf spaces: the compact, symmetric, quaternion-Kähler manifolds. Now, this construction can be extended to a larger class of inner symmetric spaces $G/K$ (see \cite{sa2}). The total space $G/H$ is a \emph{generalized flag manifold} (i.e. $H$ is the centralizer of a torus in $G$) and a 3-symmetric space. In particular $H$ contains a maximal torus of $G$ (or has maximal rank) but it is not a maximal subgroup since the inclusion $H \subset K$ is strict.

%An homogeneous space $G/H$ is called a generalized flag manifold if $H$ is the centralizer of a torus in $G$. Natural complex structures on such spaces are given by isomorphisms $G/H \simeq G^{\CM}/P$ where $P$ is a parabolic subgroup of $G^{\CM}$. To each generalized \emph{complex} flag manifold, Burstall and Rawnsley \cite{bu} associated a canonical fibration from $G^{\CM}}/P$ to a symmetric space $G/K$. This fibration is a twistor fibration, i.e. the fibres are complex submanifolds of the total space. Reciprocally,

\item[$\mathfrak B$.] The 3-symmetric spaces of the second type are those studied by Wolf \cite{wo3}, theorem 8.10.9, p280. They're characterized by $H$ being the connected centralizer of an element of order 3. Thus, $H$ is not the centralizer of a torus anymore. However it still has maximal rank, i.e. the 3-symmetric space is \emph{inner}, and is furthermore maximal (for an explicit description, using the extended Dynkin diagram of $\h$, see \cite{wo3}, \cite{wo} or \cite{bu4}).

\item[$\mathfrak C$.] Finally, the 3-symmetric spaces of the third type have ${\rm rank}\, H < {\rm rank}\, G$. Equivalently $M$ is an outer 3-symmetric space. This includes two exceptional spaces -- $Spin(8)/G_2$ and $Spin(8)/SU(3)$ -- and the infinite family $G \times G \times G/G$ of section 2.

\end{itemize}

This division, which was obtained in \cite{wo} by algebraic means (or group theory) has a profound geometrical interpretation. Indeed, the three classes can be characterized by the type of the isotropy representation :

\begin{itemize}

\item[$\mathcal A$.] In the first case, the vertical distribution $V$ and the horizontal distribution $H$ of $G/H \to G/K$ are invariant by the left action of $G$. Moreover, by definition of the natural almost complex structures, they are stable by $J_{\pm}$. Thus, the isotropy representation is \emph{complex reducible} (we identify $J_-$ with the multiplication by $i$ on the tangent spaces).

\item[$\mathcal B$.] The exceptionnal spaces that constitute the second class of 3-symmetric spaces are known to be isotropy \emph{irreducible}. Moreover, by \cite{wo3} corollary 8.13.5, these are the only non-symmetric isotropy irreducible homogeneous spaces $G/H$ such that $H$ has maximal rank.

\item[$\mathcal C$.] Finally, for the spaces of type $\mathfrak C$, the isotropy representation is reducible but not complex reducible. Let $J$ be the canonical almost complex structure of $G/H$, viewed as a constant tensor on $\m$. There exists an invariant subset $\mathfrak{n}$ such that $\m$ decomposes into $\mathfrak{n} \oplus J\mathfrak{n}$. This situation is called \emph{real reducible} by Nagy \cite{na2}.

\end{itemize}

\begin{NB}
The dimension 6 is already representative of Gray and Wolf's classification. Indeed, we have already seen that $\CM P^3$ and $\mathbb F^3$ are the twistor spaces of $S^4$ and $\CM P^2$, respectively. Secondly, the sphere $S^6 \simeq G_2/SU(3)$ is isotropy irreducible (see section 4). And thirdly, $S^3 \times S^3$ belongs to the infinite family of class $\mathfrak C$.
\end{NB}

Now, for a general nearly Kähler manifold, we can't look at the isotropy representation anymore. However we must remember proposition \ref{connexion normale = intrinseque} that the normal connection coincides with the intrinsic connection, for a 3-symmetric space and so the isotropy representation is equivalent to the \emph{holonomy representation} of $\nb$. The question becomes, then, what can we say about the geometry of a nearly Kähler manifold whose holonomy is respectively: \emph{complex reducible}, \emph{irreducible} or \emph{real reducible} ?

\begin{itemize}

\item[a.] Belgun and Moroianu, carrying out a program of Reyes Carri\'on \cite{re}, p57 (especially proposition 4.24), have shown, in \cite{be}, that the holonomy of a 6-dimensional nearly Kähler manifold $M$ is complex reducible (or equivalently the holonomy group of $\nb$ is contained in $U(2) \subset SU(3)$) if and only if $M$ is the twistor space of a positive self-dual Einstein 4-manifold. Thus, by the result of Hitchin \cite{hi}, the only compact, simply connected, complex reducible, nearly Kähler manifolds, in dimension 6 are $\CM P^3$ and $\mathbb F^3$. P.A. Nagy generalized this result in higher dimensions. Let $M$ be a complete, irreducible (in the Riemannian sense) nearly Kähler manifold of dimension $2n$, $n \geq 4$, such that the holonomy representation of $\nb$ is complex reducible. Then, $M$ is the twistor space of a quaternion Kähler manifold or of a locally symmetric space.

\item[b.] The irreducible case relates to weak holonomy. Cleyton and Swann \cite{cl,cl2} have shown an analog of Berger's theorem on special holonomies \cite{berg} for the special geometries \emph{with torsion}. By this we mean $G$-manifolds, or real manifolds of dimension $m$, equipped with a $G$-structure, $G \subset SO(m)$, such that the Levi-Civita connection of the underlying metric structure is \emph{not} a $G$-connection, or equivalently the torsion $\eta$ of the intrinsic connection is not identically zero. More precisely they made the hypothesis that $M$, or the holonomy representation of $\nb$, are irreducible and that $\eta$ is totally skew-symmetric and parallel (the last condition is automatically satisfied for a nearly Kähler manifold by proposition \ref{parallel torsion}). Then, $M$ is (\romannumeral 1) a homogeneous space or (\romannumeral 2) a manifold with weak holonomy $SU(3)$ or $G_2$. The first case in (\romannumeral 2) corresponds exactly to the 6-dimensional irreducible nearly Kähler manifolds while the second is otherwise called \emph{nearly parallel $G_2$}. Moreover, the geometry of the homogeneous spaces in (\romannumeral 1) may be specified. Indeed, the proof consists in showing that the curvature of the intrinsic connection is also parallel. Then, $\nb$ is an Ambrose-Singer connection. This reminds us of theorem \ref{loc 3symmetric}. Finally, irreducible nearly Kähler manifolds are (\romannumeral 1) 3-symmetric of type $\mathfrak B$ or (\romannumeral 2) 6-dimensional.

\item[c.] Eventually, Nagy has proved in \cite{na2}, corollary 3.1, that the complete, simply connected, real reducible, nearly Kähler manifolds are 3-symmetric of type $\mathfrak C$.

\end{itemize}

He summarized his results in a partial classification theorem (\cite{na2}, theorem 1.1): let $M$ be a complete, simply connected, (strictly) nearly Kähler manifold. Then, $M$ is an almost Hermitian product of the following spaces:
\begin{itemize}
\item 3-symmetric spaces of type $\mathfrak A$, $\mathfrak B$ or $\mathfrak C$
\item twistor spaces of non locally symmetric, quaternion-Kähler manifolds
\item 6-dimensional irreducible nearly Kähler manifolds
\end{itemize}

If we suppose furthermore that the manifold is homogeneous, then there remains only 3-symmetric spaces and 6-dimensional, nearly Kähler, \emph{homogeneous} manifolds. As a consequence, conjecture \ref{conjecture} is reduced to dimension 6.

\vs

Now, the proof of the conjecture in dimension 6 has two parts. Firstly, we must show that the only homogeneous spaces $G/H$ admitting an invariant nearly Kähler structure are those considered in sections 2, 3, 4.

\begin{lemm}
Let $(G/H,g,J)$ be a 6-dimensional almost Hermitian homogeneous space, such that the almost Hermitian structure $(g,J)$ is nearly Kähler. Then, the Lie algebras of $G$, $H$ are given at one entry of the following table. Moreover, if $G/H$ is simply connected it is isomorphic to the space at the end of the line.
\begin{equation}
\begin{array}{c|l|l|l}
\boldsymbol{{\rm dim} \, \h} & \boldsymbol{\h} & \boldsymbol{\g} & \\
\hline
0 & \{0\} & \mathfrak{su}(2) \oplus \mathfrak{su}(2) & S^3 \times S^3 \\
\hline
1 & i\RM & i\RM \oplus \mathfrak{su}(2) \oplus \mathfrak{su}(2) & S^3 \times S^3 \\
\hline
2 & i\RM \oplus i\RM & i\RM \oplus i\RM \oplus \mathfrak{su}(2) \oplus \mathfrak{su}(2) & S^3 \times S^3 \\
  & i\RM \oplus i\RM & \mathfrak{su}(3) & \mathbb F^3 \\ 
\hline
3 & \mathfrak{su}(2) & \mathfrak{su}(2) \oplus \mathfrak{su}(2) \oplus \mathfrak{su}(2) & S^3 \times S^3 \\
\hline
4 & i\RM \oplus \mathfrak{su}(2) & i\RM \oplus \mathfrak{su}(2) \oplus \mathfrak{su}(2) \oplus
\mathfrak{su}(2) & S^3 \times S^3 \\
  & i\RM \oplus \mathfrak{su}(2) & \mathfrak{sp}(2) & \CM P^3 \\
\hline
8 & \mathfrak{su}(3) & \g_2 & S^6 \\
\hline
\end{array}
\label{table h,g}
\end{equation}
\label{algebraic classification}
\end{lemm}

\begin{proof}
By a result of Nagy \cite{na}, the Ricci tensor of a nearly Kähler manifold $M$ is positive (in dimension 6, this a consequence of Gray's theorem in \cite{gr3} that $M$ is Einstein, with positive scalar curvature). Then, by Myer's theorem, $M$ is compact with finite fundamental group and the universal cover $\widetilde M$ of $M$ is a nearly Kähler manifold of the same dimension. Moreover, if $M \simeq G/H$ is homogeneous, $\widetilde M$ is isomorphic to $\widetilde G/ \widetilde H$ where the groups $G$, $\widetilde G$ and $H$, $\widetilde H$ have the same Lie algebras. Consequently, we shall work with $\widetilde M$ instead of $M$.

For an homogeneous space, we have the following homotopy sequence:
\begin{equation}
\dots \to \pi_2(G/H) \to \pi_1(H) \to \pi_1(G) \to \pi_1(G/H) \to H/H^0 \to 0
\label{homotopie}
\end{equation}
If the manifold is simply connected, this provides us for a surjective automorphism $\varphi$, from the fundamental group of $H$ to the fundamental group of $G$. We shall look only at the consequences at the Lie algebra level (i.e. at the $S^1$ or $i\RM$ factors and not at the finite quotients).

The second feature we use is the natural reduction to $SU(3)$ defined in section 2. If $g$ and $J$ or $\omega$ are invariant for the left action of $G$, then, so is $\psi = d\omega$. As a consequence, the isotropy group $H$ is a subgroup of $SU(3)$. This leaves the following possibilities for $\h$: $\{0 \}$, $\mathfrak{u}(1)$, $2\mathfrak{u}(1)$, $\mathfrak{su}(2)$, $\mathfrak{u}(2)=\mathfrak{su}(2) \oplus \mathfrak{u}(1)$ and $\mathfrak{su}(3)$. Next, to find $\g$, we use the fact that the difference between the dimensions of the two algebras is $6$, the dimension of the manifold, and the existence of $\varphi$ above. The latter allows us to eliminate the following : $\g = \mathfrak{su}(2) \oplus 3\mathfrak{u}(1)$ or $6\mathfrak{u}(1)$ for $\h = \{0\}$ and $\g=\mathfrak{su}(3) \oplus \mathfrak{u}(1)$ for $\h = \mathfrak{su}(2)$. Table (\ref{table h,g}) is a list of the remaining cases.

Now, $\h$ acts as a subgroup of $\mathfrak{su}(3)$ on the 6-dimensional space $\m$. Using this, we determine the isotropy representation and the embedding of $H$ into $G$. In particular we show, when $\g = \h \oplus \mathfrak{su}(2) \oplus \mathfrak{su}(2)$, that $G/H$ is isomorphic to $S^3 \times S^3$ and $G$ contains $SU(2) \times SU(2)$ acting on the left. So the nearly Kähler structures arising from these cases will always induce a left-invariant nearly Kähler structure on $S^3 \times S^3 \simeq SU(2) \times SU(2)$. Thus, we need only to consider this more general situation. This was done in section 2.
\end{proof}

Now, theorem \ref{classification} is a consequence of lemma \ref{algebraic classification} and propositions \ref{S3xS3}, \ref{CP3}, \ref{F3} and corollary \ref{S6}.

\vs

\small\noindent Université de Cergy-Pontoise \\
Département de Mathématiques\\
site de Saint-Martin\\
2, avenue Adolphe Chauvin\\
95302 Cergy-Pontoise cedex, France\\
\texttt{Jean-Baptiste.Butruille@u-cergy.fr}

\end{document}